\documentclass[12pt]{amsart}
\usepackage{latexsym,amsfonts,amsmath,amssymb,amsthm,rotating}
\usepackage{epic}
\newtheorem{theorem}{Theorem}[section]
\newtheorem{cor}[theorem]{Corollary}

\newtheorem{prop}[theorem]{Proposition}
\theoremstyle{definition}

\newcommand{\bm}{\mathbf{m}}

\newcommand{\C}{{\mathbb C}} 
\newcommand{\Z}{{\mathbb Z}} 
 
\newcommand{\even}{{\mathrm{even}}}
\newcommand{\odd}{{\mathrm{odd}}}

\newcommand{\ptra}{\mathrm{PT}}
\newcommand{\oar}[1]{\overset{#1}{\longrightarrow}}
\newcommand{\oal}[1]{\overset{#1}{\longleftarrow}}
\newcommand{\rel}{\mathrm{Rel}}
\newcommand{\Lambdat}{\hat\Lambda}
\newcommand{\Sigmat}{\widehat\Sigma}
\newcommand{\ot}[1]{_{\oplus#1}}
\newcommand{\sigmat}{\hat\sigma}
\title{Periodicity of Y-systems and flat connections} 
\author{ Andr\'as Szenes} \date{}

\begin{document}
\maketitle

\section{Introduction}

The $Y$-system is a multi-dimensional rational recursion which first
appeared in theoretical physics. It was introduced by Al.
Zamolodchikov \cite{zam}, who was motivated by ideas coming from
deformed conformal field theories (see also \cite{kuna} for
connections with integrable systems). The recursion was generalized by
Ravanini, Tateo and Valleriani in \cite{Ravanini}, and in this general
form it is parameterized by a pair of Dynkin diagrams.  In the present
paper we restrict our attention to the case of Dynkin diagrams of type
$A$.

The recursion starts with a certain number of free parameters and
progresses along according a set of relations (see \eqref{zam_eq}).
The conjecture is that the recursion eventually returns to its
original starting conditions (Theorem \ref{thethm}). To prove this
turned out to be surprisingly difficult, in part, because it is not
quite clear what mathematical tools one could employ. The case
$A_1\times A_k$ was proved by Frenkel-Szenes \cite{fsz} and
Gliozzi-Tateo \cite{gt}, using a rational parameterization based on
continued fractions and hyperbolic geometry, respectively.

The cases $A_1\times \Delta$, where $\Delta$ is any Dynkin diagram has
been proved by Fomin and Zelevinsky \cite{fz}, where they successfully
related the system to their theory of cluster algebras. This case
provided one of the basic examples of the theory. In a remarkable
recent paper \cite{fzrec}, Fomin and Zelevinsky also show that
replacing the exponents in the recursion in the $A_1\times M$ case
with a matrix $M$ in a rather wide class, the system will be periodic
only if $M$ comes from a Dynkin diagram.

We also note that the $Y$-system of Zamolodchikov has been linked to
identities of the dilogarithm functions  (\cite{gt},\cite{fsz}
\cite{chap}).

The general case of the product of two Dynkin diagrams has been open
for more 10 years now, even in the case $A_k\times A_r$. In this case
the $Y$-system is not clearly related to the Fomin-Zelevinsky theory
cluster algebras. 

In this paper, we give a proof of the periodicity of the $A_k\times
A_r$ system, using a novel interpretation of the system as a system of
flat connections on a graph.

When this work substantially completed we learned of another proof of
the periodicity for this case by Alexandre Volkov \cite{volk}.
Volkov's proof is rather different; it uses an explicit
parameterization.  The two proofs are of independent interest. It
would be interesting to see if either of these proofs could be
generalized to the arbitrary Dynkin diagram case.

We should mention that a similar periodicity phenomenon has been
observed by A. Henriques \cite{enr}.

{\bf Acknowledgments}. We are greatly indebted to the ideas and generous
help of A. Volkov. We would also like to express our gratitude to the
University of Geneva for their hospitality, and to Anton
Alekseev for his advice and encouragement.

\section{The infinite system}

Consider the 3-dimensional lattice $\Lambda=\Z\times\Z\times\Z$. We
will need to visualize this lattice, and our convention will be that
the first coordinate is horizontal (East-West), the second is vertical
(North-South, increasing towards the South), and the third is up-down
(plus-minus). The infinite $Y$-system is an algebraic system with
infinitely many variables and infinitely many relations, both indexed
by the lattice points in $\Lambda$.

It is best to think of the variables as formal ones, but for
simplicity of notation we will consider them as having actual
complex values. Thus consider the space of complex valued functions on $\Lambda$
\[ F(\Lambda) = \{Y:\Lambda\to\C^*\}.\]
Define the shift operators on this space: for $\bm=(n,i,j)\in\Lambda$, let
\begin{multline*}
Y^W(\bm)=Y(n-1,i,j),\,Y^N(\bm)=Y(n,i-1,j),\,Y^E(\bm)=Y(n+1,i,j),\\
Y^S(\bm)=Y(n,i+1,j),\,Y_+(\bm)=Y(n,i,j+1),\,Y_-(\bm)=Y(n,i,j-1).
\end{multline*}

Introduce the so-called $Y$-{\em system}:
\begin{equation}
  \label{zam_eq}
Y^WY^E=\frac{(1+Y^{N})(1+Y^{S})}{(1+1/Y_+)(1+1/Y_-)}.  
\end{equation}
Thus indeed, for each site $\bm\in\Lambda$, there is one value,
$Y(\bm)$, and also one relation, $\rel(\bm)$, relating
the values of $Y$ on the 6 neighboring sites.

Three quick observations about this system.
\\ {\bf Decoupling}: The lattice $\Lambda$ is the disjoint union of
the lattices
\[ \Lambda^{\even}=\{(n,i,j);\;i+j+n\text{ even}\}\quad\text{and}\quad
\Lambda^{\odd}=\{(n,i,j);\;i+j+n\text{ odd}\}.
\]
The system clearly decouples into two, one with variables
parameterized by $\Lambda^{\even}$:
\[   F^{\even}=\{Y:\Lambda^\even\to\C^*\},
\]
and relations parameterized by $\Lambda^\odd$, and another, equivalent
system with the roles of {\em even} and {\em odd} exchanged. From now
on we will only consider the even system.  \\ {\bf Rationality}:
Again, by inspecting the relations \eqref{zam_eq}, one can easily see
that the variables $\{Y(0,i,j),Y(1,i,j);\;i,j\in\Z\}$ are independent
and determine the rest of the variables uniquely. More precisely,
using \eqref{zam_eq}, any variable of the system may be expressed as a
rational function of the variables with $n=0,1$. Moreover, a generic
function $\hat Y:\{(n,i,j);\;n=0\text{ or }1\}\to\C^*$ extends to a
unique element of $F(\Lambda)$ satisfying \eqref{zam_eq}.
\\
{\bf Symmetry}: By replacing $Y$ with $1/Y$ the system turns into
another one, which has exactly the same form, with the North-South
direction exchanged with the up-down direction.

\section{The truncation and the formulation of the conjecture}
Now we introduce a truncated version of \eqref{zam_eq}. Let
\[\Lambda_{rk}=\{(n,i,j);\;1\leq i\leq r,\, 1\leq j \leq k\},
\] be the truncated lattice and denote  by $F(\Lambda_{rk})$ the
corresponding function space. For each $\bm\in\Lambda_{rk}$ we impose
the relation $\rel(\bm)$ given in \eqref{zam_eq}, with the {\em
  convention} that the factors with $i,j<1,\,i>r$ and $j>k$, are
simply omitted from the equations. Another way to put this is to set the
boundary conditions
\[ Y(n,0,j)=Y(n,r+1,j)=1/Y(n,i,0)=1/Y(n,i,k+1)=0.
\]

Just as in the infinite case, it is not hard to see that this system
still has the decoupling, rationality and symmetry properties, where
this latter one now exchanges $r$ and $k$. In particular, we
will only consider the even system $F(\Lambda_{rk}^\even)$, whose
relations are parameterized by $\Lambda^\odd_{rk}$, and without loss
of generality we can assume that $r\geq k$.

Now we are ready to formulate the periodicity conjecture.
Define a twisted shift in the truncated lattice $\Lambda_{rk}$ via the
formula
\[ \sigma:(n,i,j)\mapsto(n+r+k+2,r+1-i,k+1-j).\]
This is thus a combination of two central symmetries of the Dynkin
diagrams $A_r$ and $A_k$, and a translation in the $\Z$-direction by
$r+k+2$.  

\begin{theorem} \label{thethm}
  Consider a solution of the truncated $Y$-system, which is a map
  $Y:\Lambda_{rk}\to\C^*$ satisfying the equations \eqref{zam_eq} with
  the convention described above. Then we have $Y^\sigma=Y$, i.e.
\[
Y(n,i,j)=Y(n+r+k+2,r+1-i,k+1-j).
\]
\end{theorem}
Note that $\sigma^2$ is simply translation by $2(r+k+2)$, hence the
system is indeed periodic.

Recall that according to the rationality property above, starting from
a generic set of values of $[Y(0,i,j),Y(1,i,j)]$ we obtain a solution
of the truncated $Y$-system, but it is not at all clear why such
system needs to be periodic.

\section{Flat connections}
\subsection{The shifted system} It is convenient to introduce the following
shifted form of our system. Let $z_{j}(n,i+j)=-Y(n,i,j)$. Our system
\eqref{zam_eq} now looks as follows:
\begin{equation}
  \label{zsys}
z^Wz^E = \frac{1-z^N_-}{1-1/z^N}\frac{1-z^S_+}{1-1/z^S}. 
\end{equation}

In this paragraph, we carefully and mechanically, rewrite all the
statements and observations of the previous section in this shifted
form. While this may seem tedious, without these formulas the proof
will be difficult to follow.

To define the property of {\bf decoupling}, consider the lattices
\[ \Lambdat^{\even}=\{(n,i,j);\;i+n\text{ even}\}\quad\text{and}\quad
\Lambdat^{\odd}=\{(n,i,j);\;i+n\text{ odd}\}.
\]
Then again the system decouples into two, with the variables sitting
on $\Lambdat^\even$ related by equations indexed by $\Lambdat^\odd$.

The truncated configurations look more complicated. We have the
truncated lattices
\[
\Lambdat^{\even}_{rk}=\{(n,i,j)\in\Lambdat^\even;\;1\leq j\leq
k,\;j+1\leq i\leq j+r\}
\]
and
\[
\Lambdat^{\odd}_{rk}=\{(n,i,j)\in\Lambdat^\odd;\;1\leq j\leq
k,\;j+1\leq i\leq j+r\}.
\]

The index $j$ thus takes the values $j=1,\dots,k$, and for $j=1$ we
have
\begin{equation}
  \label{jtrunc0}
z^Wz^E = \frac{1}{1-1/z^N}\frac{1-z^S_+}{1-1/z^S},    
\end{equation}
while for $j=k$,
\begin{equation}
  \label{jtrunck}
z^Wz^E = \frac{1-z^N_-}{1-1/z^N}\frac{1}{1-1/z^S}.
\end{equation}
The situation for the index $i$ is a bit different: it takes the
values $i=j+1,\dots,j+r$, and for $i=j+1$ we have
\begin{equation}
  \label{itrunk0}
z^Wz^E = (1-z^N_-)\frac{1-z^S_+}{1-1/z^S},  
\end{equation}
while for $i=j+r$
\begin{equation}
  \label{itruncj}
z^Wz^E = \frac{1-z^N_-}{1-1/z^N}(1-z^S_+),  
\end{equation}

The 4 additional cases for the 4 edges are similar, with two factors
omitted. For example, for $j=1,\,i=j+1$, we have
\begin{equation}
  \label{edgeeq}
z^Wz^E = \frac{1-z^S_+}{1-1/z^S}.
\end{equation}

\subsection{The $\Gamma$-system} 

An important step in our proof is the reinterpretation of the
$z$-system as the conditions of flatness of a connection defined on a
certain graph. Again, first we consider the simpler infinite case.

We introduce another infinite system now, which we will call the
$\Gamma$-system. Formally, the variables of this system live on the
points of the shifted lattice $v+\Lambda$, where $v=(-1/2,-1/2,0)$.
We will use different symbols and conventions for the variables
corresponding to $v+\Lambdat^\even$ and $v+\Lambdat^\odd$. Assuming
that $n+i$ is even, we will write 
\[a_j(n,i)\text{ for the variable at }(n+1/2,i+1/2,j)\in
v+\Lambdat^\even
\]
and
\[x_j(n,i)\text{ for the variable at }(n-1/2,i+1/2,j)\in
v+\Lambdat^\odd
\]
We can also think of each variable of this new system as sitting on an
interval joining two points of $\Lambdat^\odd$; in this picture then
the points of $\Lambdat^\even$ are the centers of squares formed by
these intervals.

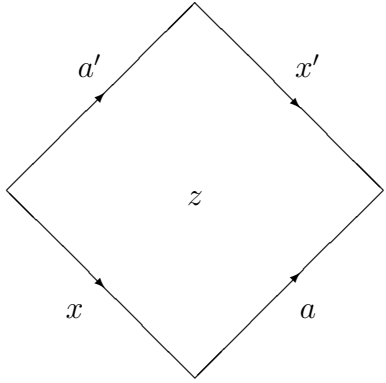
\begin{center}

\setlength{\unitlength}{1mm}
\begin{picture}(60,60)(-5,-5)
\linethickness{1pt}
  \put(0,25){\line(1,1){25}}
\put(25,50){\line(1,-1){25}}
\put(50,25){\line(-1,-1){25}}
\put(25,0){\line(-1,1){25}}

\put(40,8){\makebox(0,0)[b]{$a$}}
\put(11,40){\makebox(0,0)[b]{$a'$}}
\put(9,8){\makebox(0,0)[b]{$x$}}
\put(40,40){\makebox(0,0)[b]{$x'$}}
\put(25,23){\makebox(0,0)[b]{$z$}}

   \put(10,35){\vector(1,1){3}}
\put(36,39){\vector(1,-1){3}}
\put(36,11){\vector(1,1){3}}
\put(10,15){\vector(1,-1){3}}

  \end{picture}
  \end{center}

Now consider the directed graph $\Gamma$ obtained as the projection of
these intervals onto the $ni$-plane, with directions chosen eastward:
northeast for the $a$-edges, and southeast for the $x$-edges.
Associate to each SE edge the infinite matrix $X(i,j)$ which has 1s on
the diagonal and $x_j(n,i)$ as the $(j,j+1)$ entry. Also, associate to
NE edge the matrix $A(n,i)$, which has $a_j(n,i)$ as the $j$th
diagonal entry and all the entries under the diagonal are 1s. The
rest of the entries vanish. 
These matrices define a connection on the graph $\Gamma$ with values
in $GL(\infty)$, and the {\em equations of the $\Gamma$-system} say that
this connection is flat, i.e.
\begin{equation}
  \label{gammasys}
XA = A'X',  
\end{equation}
where we used the notation $A'=A^{NW}$ and $X'=X^{NE}$.

Explicitly, in terms of the entries, the equations look like
\begin{equation}
  \label{xarels}
x+a = a'+x_-',\quad xa_+=a'x'.  
\end{equation}
Here we think of $x$ as a function on $\Lambdat^\even$, $a$ is
a function on $\Lambdat^\odd$, and again, $a'=a^{NW}$ and $x'=x^{NE}$.

Again, we have the rationality property of this system: now this means
that the values of $x(0,i,j)$ and $a(0,i,j)$ are independent, and that
all other values of $x$ and $a$ are rational functions of these values.

\begin{prop}\label{gammatoz}
  Any $\Gamma$-system gives rise to a $z$-system via
  the formula
  \begin{equation}
    \label{trafo}
z = \frac x{a'},    
  \end{equation}
and all $z$-systems arise this way. 
\end{prop}

  \begin{center}
  \setlength{\unitlength}{1mm}
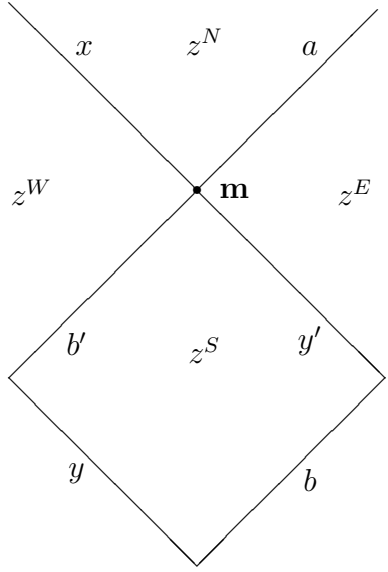
\begin{picture}(50,80)(0,-25)
\linethickness{1pt}
  \put(0,0){\line(1,1){49}}
\put(0,50){\line(1,-1){50}}
 \put(0,0){\line(1,-1){25}}
\put(50,0){\line(-1,-1){25}}
  \put(25,25){\circle*{1}}
  \put(40,43){\makebox(0,0)[b]{$a$}}
\put(10,43){\makebox(0,0)[b]{$x$}}
\put(9,3){\makebox(0,0)[b]{$b'$}}
\put(40,3){\makebox(0,0)[b]{$y'$}}
\put(9,-14){\makebox(0,0)[b]{$y$}}
\put(40,-15){\makebox(0,0)[b]{$b$}}

\put(26,43){\makebox(0,0)[b]{$z^{N}$}}
\put(26,2){\makebox(0,0)[b]{$z^{S}$}}  
\put(46,23){\makebox(0,0)[b]{$z^{E}$}}
\put(30,24){\makebox(0,0)[b]{$\bm$}}
\put(3,23){\makebox(0,0)[b]{$z^{W}$}}
  \end{picture} 
  \end{center}

\begin{proof}
Let $(x,a)$ be a solution of the $\Gamma$-system, and let $\bm$ be a
lattice point in $\Lambdat^\odd$ marking one of the relations in
\eqref{zsys}. To simplify our notation, we will denote the variables
surrounding $z^N$ by $x,a,x',a'$ as usual, while for the variables
surrounding $z_S$ we will use $y,b,y',b'$.

Now we compute the factors in \eqref{zsys}. We have
\[  \frac{1-z_+^S}{1-1/z^S}=\frac{1-y_+/b_+'}{1-b'/y}=
\frac{b_+'-y_+}{y-b'}\cdot\frac{y}{b_+'}=-\frac{y}{b_+'}.
\]
where we used the first equation in \eqref{xarels}.
For the other case, observe that the multiplicative equation in
\eqref{xarels} maybe rewritten as $x/a'=x'/a_+$. Then
\begin{equation}
  \label{xpera}
\frac{1-z^N_-}{1-1/z^N}=\frac{1-x_-'/a}{1-a'/x}=\frac{a-x_-'}{x-a'}\cdot
 \frac xa = -\frac xa  
\end{equation}
Finally, note that 
\[  z^W =\frac x{b_+'},\quad\text{and}\quad z^E=\frac{y'}a.
\]
Now the equation \eqref{zsys} immediately follows.
\end{proof}

Thus we have managed to interpret the $Y$-system as a system flat
graph connections. Our next goal is to  impose boundary
conditions on the $\Gamma$-system in such a way that they induce the
truncation of the $z$-system. 

Imposing the condition in the $j$-direction is very natural: consider
the $\Gamma$-system \eqref{gammasys}, with $k+1$-by-$k+1$ matrices
instead of infinite ones. Thus now we have
\[ X(n,i) = 
\begin{pmatrix}
  1 & x_1 & 0 & \hdotsfor{2}&  0 \\
  0 & 1 & x_2 & 0 & \dots & 0 \\
\hdotsfor{6}\\
0 & \dots & 0 & 1 & x_{k-1}& 0 \\
0 & \hdotsfor{2} & 0 & 1 & x_{k}\\
0 & \hdotsfor{3} & 0 & 1 
  \end{pmatrix}
\]
and
\[ A(n,i) = 
\begin{pmatrix}
  a_1 & 0 & 0 & \hdotsfor{3}&  0 \\
   1 & a_2 & 0 & \hdotsfor{3} & 0 \\
\hdotsfor{7}\\
0 & \dots & 0 & 1 &  a_{k-1}& 0 & 0 \\
0 & \hdotsfor{2} & 0 & 1 & a_{k}& 0\\
0 & \hdotsfor{3} & 0 & 1 & a_{k+1} 
  \end{pmatrix}
\]

For this truncated system the index $j$ for the $a$ variable runs
from $j=1$ to $k+1$, while for the variable $x$ this interval is $j=1,\dots,k$.

This truncation of the $\Gamma$-system means that for $j=1$ we have
$x+a=a'$, while for $j=k+1$ we have $a=a'+x_-'$. Now it is easy to
check that the two equations \eqref{jtrunc0} and \eqref{jtrunck} follow
from this truncated $\Gamma$-system exactly as in the proof of
Proposition \ref{gammatoz}

Now we come to the truncation in the $i$ direction. Here the situation
will be a not as pleasant, because we will not be able to interpret
the resulting system as a flat connection on a graph immediately. 

The truncated $\Gamma$-system is defined as follows:
\[ x_j(n,i)\text{ is defined for } 1\leq j \leq k,\; j\leq i\leq j+r,
\]
while
\[ a_j(n,i)\text{ is defined for }
\begin{cases}
   j=1,\; 1\leq i \leq r+1\\
  1\leq j\leq k,\;j\leq i \leq j+r-1\\
j=k+1,\; k\leq i \leq k+r
\end{cases}
\]
The relations \eqref{xarels} are truncated as follows:
\begin{equation}
  \label{truncgam}
\begin{cases}
  x+a=a'\text{ if }j=1,\\
  a=a'+x_-'\text{ if }j=k+1,\\
  x+a=x_-'\text{ if }i=j,\\
  x=a'+x_-'\text{ if }i=j+r.
\end{cases}  
\end{equation}

Now we claim
\begin{prop}
  Consider the truncated $\Gamma$-system described above. Then the
  system $z=x/a'$ satisfies the equations of the truncated $z$-system.
  The map from truncated $\Gamma$-systems to the truncated $z$-systems
  is surjective.
\end{prop}

The proof is identical to that of Proposition \ref{gammatoz}. We would
like to note that the variables $x_j(0,i)$, $a_j(0,i)$ are still
independent, but now they do not necessarily determine the rest of the
variables. Sometimes one has the freedom of choosing some of the edge
$a$ variables as the recursion progresses, but this fact does not
influence our result.

\section{Boundary conditions and permutations}

As we mentioned above, the problem is that these truncated
$\Gamma$-systems are cannot be interpreted as flat connections since
for certain values of $i$ not all entries of the corresponding
matrices are defined.

The key idea is that our boundary conditions at these partially
defined matrices force a certain matrix transformation which then
quickly leads to the proof of the periodicity of the system.

We demonstrate this transformation for the $k=2$ case first. The
general case is analogous. It is described in Proposition \ref{propeq}.

Thus we have a truncated system with $k=2$, i.e. with 3-by-3 matrices,
and consider the $i=1$ case, which is the critical one, because the
system is not defined for $j=2$. In other words, fixing $i=2$, we have
$x_1,x_2, a_1,a_2,a_3$ and $a_1',x_1'$ but not $a_2',x_2',a_3'$. Now
we are unable to write down our usual matrix equation $XA=A'X'$.
However, taking into account the truncated $\Gamma$-equations
\eqref{truncgam}, concretely
\[ x_1+a_1=a_1',\; x_2+a_2=x_1',\; x_1a_2=a_1'x_1',
\]
we have
\[ XA = 
\begin{pmatrix}
  a_1' & 0 & 0 \\
1 & x_2 & 0 \\
0 & 0 & 1 
\end{pmatrix}
\begin{pmatrix}
  1 & 0 & 0 \\
0 & 0 & a_3 \\
0 & 1 & a_3
\end{pmatrix}
\begin{pmatrix}
  1 & x_1' & 0 \\
0 & 1 & 0 \\
0 & 0 & 1
\end{pmatrix}
\]
This equality can be checked by multiplying the matrices and comparing
the corresponding entries. Clearly, this equality is local, i.e. we
can write it down for any $k$. The structure is described in the
following statement.

\begin{prop}\label{propeq}
For $i=k$ we will have 
\begin{equation}
  \label{geneq}
XA=A'\ot1 \Sigmat_{k,k+1}X'\ot1  
\end{equation}
where $A'$ and $X'$ are $k$-by-$k$ matrices, whose entries (apart from
the lower-left entry of $A'$) are exactly those coming from the $i=k-1$ part
of the system, and the matrix $\Sigmat_{k,k+1}$ is identity on the
first $k-1$ coordinates and a 2-by-2 transformation in the last two
coordinates with the $(k,k)$ entry vanishing.
\end{prop}
Here we used the notation $M\ot d$ taking the direct sum of a matrix
$M$ with the identity matrix in dimension $d$. When $d=1$, this means
adding a row and a column of zeros to $M$ and then changing the
lower-right entry to 1.

By iterating this transformation, we arrive at the following key fact:
\begin{prop}\label{smatrix} Let $(x,a)$ be a truncated
  $\Gamma$-system, pick any $n$ such that $n+k$ is even. Then  the product
  of matrices along the staircase path 
  \begin{equation}
    \label{smat}
X(n,k)A(n,k)X(n+2,k)A(n+2,k)\dots X(n+2(k-1),k)A(n+2(k-2),k)    
  \end{equation}
  is anti-lower-triangular, i.e. it has the property that all of its
  $[j,l]$ entries with $j+l<2k$ vanish.
\end{prop}
Remark: A similar statement holds for the Southern edge: the product
\[ A(n,r+1)X(n,r+1)\dots X(n+2(k-1),r+1)A(n+2(k-2),r+1)
\]
is anti-upper-triangular.

\begin{proof}
Let us demonstrate the proof for the case $k=3$. Without loss of
generality we can choose $n=-1$. Then we start with
\[ X(-1,3)A(-1,3)X(1,3)A(1,3) X(3,3)A(3,3)    
\]
and apply the \eqref{geneq} to each of the three  products of the form
$XA$. We end up with
\[ A(-2,2)\ot1\Sigmat_{34}X(0,2)\ot1
A(0,2)\ot1\Sigmat_{34}X(2,2)\ot1
A(2,2)\ot1\Sigmat_{34}X(4,2)\ot1,
\]
where somewhat loosely, we denoted all the matrices $\Sigmat_{34}$ by
the same symbol, even though only their shapes coincide. Observe that
an anti-lower-triangular matrix will remain such after multiplication
by $A$ on the left or by $X$ on the right, thus we can omit these
matrices from the formula.  Now we proceed to apply the equality
\eqref{geneq} for $k=2$ to the two products of the form $XA$ in this
formula. We obtain
\[ \Sigmat_{34}
A(-1,1)\ot2\Sigmat_{23}X(1,1)\ot2
\Sigmat_{34}
A(1,1)\ot2\Sigmat_{23}X(3,1)\ot2
\Sigmat_{34}
\]
Observe that because of their direct sum structure the matrices
$\Sigmat_{34}$ commute with any matrix of the form $M\ot2$. Finally,
noting that for $k=1$ the equality \eqref{geneq} simply means that
$X(1,1)\ot2A(1,1)\ot2=\Sigmat_{12}$ we can rewrite \eqref{smat} as
\[ L \Sigmat_{34}
\Sigmat_{23}\Sigmat_{12}
\Sigmat_{34}
\Sigmat_{23}
\Sigmat_{34} U,
\]
where $L$ is a lower-, $U$ is an upper-triangular matrix.
This formula mimics the standard way of writing the longest element
of the symmetric group as a product of nearby transpositions. We leave
it as an exercise to check that such a matrix is
anti-lower-triangular.
  \end{proof}

  This result has the following beautiful corollary.  Because of the
  symmetry property of the system, we can assume
  without loss of generality that $r\geq k$. Recall that in our
  truncated system, the values of $i$ will vary between $1$ and
  $k+r+2$, however, those sites for which the matrices $X(n,i)$ and
  $A(n,i)$ are defined have $k\leq i \leq r+1$.  We will call the
  edges of $\Gamma$ for which the matrices are defined {\em regular}.
  Also, we will call {\em regular} those vertices of $\Gamma$ which
  are endpoints of these edges.

Just as in the infinite case, we can define the parallel transport
$\ptra(p_1,p_2)$ between two regular vertices $p_1,p_2$, even when the
$\Gamma$-system is truncated, by multiplying the $k+1$-by-$k+1$
matrices from one point to the other along any path consisting of
regular edges.  Again, this is well-defined because of the flatness
condition.

Now we can formulate an important corollary of Proposition \ref{smatrix}.
\begin{cor}\label{maincor}
  Define a transformation of the $ni$-plane via the formula
\[ \sigmat(n,i)=(r+k+2+n,r+k+2-i).
\]
Pick a regular vertex $p$ of $\Gamma$; then the vertex $\sigmat p$ is
also regular, and the parallel transport $\ptra(p,\sigmat p)$ is
anti-diagonal.
\end{cor}

The following picture, showing the $k=2$, $r=3$ case, will be helpful.
Here the regular part of $\Gamma$ is shown with thin solid lines

  \begin{center}
  \setlength{\unitlength}{1mm}
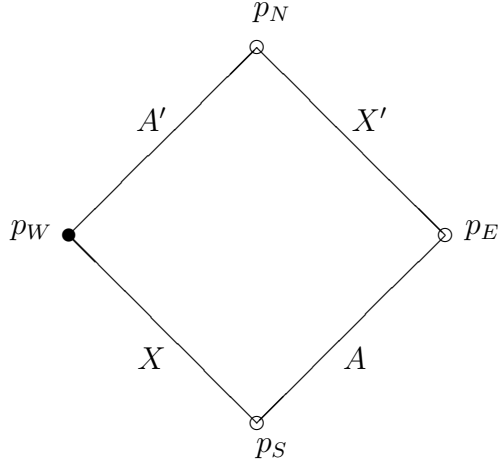
\begin{picture}(80,70)(-5,-5)
\linethickness{.8pt}
 
 \put(0,60){\line(1,0){20}}
 \put(0,50){\line(1,0){30}}
\put(10,40){\line(1,0){30}}
\put(20,30){\line(1,0){30}}
\put(30,20){\line(1,0){30}}
\put(40,10){\line(1,0){20}}

 \put(0,60){\line(0,-1){10}}
 \put(10,60){\line(0,-1){20}}
\put(20,60){\line(0,-1){30}}
\put(30,50){\line(0,-1){30}}
\put(40,40){\line(0,-1){30}}
\put(50,30){\line(0,-1){20}}
\put(60,20){\line(0,-1){10}}

\dashline{2}(0,50)(0,40)(10,40)(10,30)(20,30)(20,20)(30,20)(30,10)(40,10)

\dashline{2}(20,60)(30,60)(30,50)(40,50)(40,40)(50,40)(50,30)(60,30)(60,20)

\linethickness{1.5pt}
\put(10,50){\line(0,-1){10}}
\put(10,50){\line(1,0){20}}
\put(50,20){\line(0,1){10}}
\put(50,20){\line(-1,0){20}}

\put(30,40){\line(1,0){10}}
\put(40,30){\line(1,0){10}}
\put(10,40){\line(1,0){10}}
\put(20,30){\line(1,0){10}}

\put(30,20){\line(0,1){10}}
\put(20,30){\line(0,1){10}}
\put(40,30){\line(0,1){10}}
\put(30,40){\line(0,1){10}}


 \put(10,50){\circle*{1.8}}
   \put(7,51){\makebox(0,0)[b]{$p$}}
 \put(3,21){\makebox(0,0)[b]{{\large S}}}

  \put(50,20){\circle*{1.8}}
   
       \put(53,15){\makebox(0,0)[b]{$\sigmat p$}}
       \put(57,45){\makebox(0,0)[b]{{\large N}}}

\end{picture}
\end{center}

Introduce the notation $\delta(p) = \ptra(p,\sigmat p)$.
Consider the two paths between $p$ and $\sigmat p$ shown on the
picture. Taking the northern path, we see that $\delta(p)$ has a
representation
\[ \delta(p) =  A_1\dots A_{l}\bar LX_1\dots X_{l'} 
\]
where $\bar L$ is is a matrix of the form \eqref{smat}, thus it is
anti-lower-triangular. This implies that $\delta(p)$ is also
anti-lower-triangular.

On the other hand, going along the southern path, we have
\[ \delta(p) =   Y_1\dots Y_{l'}\bar UB_1\dots B_{l}
\]
where $\bar U$ is  anti-upper-triangular; this implies that
$\delta(p)$ is also anti-upper-triangular.

The two statements together imply that $\delta(p)$ is anti-diagonal. \qed

\section{The completion of the proof}

First we reformulate the periodicity property for the shifted system.
Let us extend the transformation $\sigmat$ from the $ni$-plane to the
whole lattice $\Lambda$ via
\[ \sigmat(n,i,j)=(r+k+2+n,r+k+2-i,k+1-j).
\]
Then the periodicity property is equivalent to saying that the
$z$-system is invariant under this transformation.

We will show this by proving an invariance property for the
corresponding truncated $\Gamma$-system.

We again assume that $r\geq k$. Consider a square in the graph
$\Gamma$ consisting of regular edges. Introduce the following notation
for the vertices and the edge matrices:
  \begin{center}
  \setlength{\unitlength}{1mm}
\begin{picture}(60,65)(-5,-8)
\linethickness{1pt}
\put(0,25){\line(1,1){25}}
\put(25,50){\line(1,-1){25}}
\put(50,25){\line(-1,-1){25}}
\put(25,0){\line(-1,1){25}}

 \put(0,25){\circle*{1.8}}
 \put(25,50){\circle{1.8}}
 \put(50,25){\circle{1.8}}
 \put(25,0){\circle{1.8}}

 \put(38,10){\makebox(0,0)[t]{$A$}}
\put(11,10){\makebox(0,0)[t]{$X$}}
\put(11,42){\makebox(0,0)[t]{$A'$}}
\put(40,42){\makebox(0,0)[t]{$X'$}}

\put(27,56){\makebox(0,0)[t]{$p_{N}$}}
\put(27,-2){\makebox(0,0)[t]{$p_{S}$}}  
\put(55,27){\makebox(0,0)[t]{$p_{E}$}}
\put(-5,27){\makebox(0,0)[t]{$p_{W}$}}

  \end{picture}
\end{center}

Thus we have
\[p_W\oar{A'} p_N\oar{X'}p_E\oal Ap_S\oal Xp_W.\]

Similarly, we set the notation
\[\sigmat p_W\oar{Y}\sigmat p_N\oar{B}\sigmat p_E\oal{B'}\sigmat p_S\oal{A'}
\sigmat p_W.\]

We have the obvious equalities
\begin{eqnarray} \label{eqn4}
& A'\delta(p_N)=\delta(p_W)Y\\ \nonumber
& A\delta(p_E)=\delta(p_S)Y' \\ \nonumber
& X\delta(p_S)=\delta(p_W)B' \\
& X'\delta(p_E)=\delta(p_N)B. \nonumber
\end{eqnarray}
Denoting the diagonal part of a matrix $M$ by $\Delta{M}$, and its
anti-diagonal part by $\bar\Delta M$, we see that
\begin{eqnarray} \label{eqn5}
& \Delta A'\bar\Delta\delta(p_N)=\bar\Delta\delta(p_W)\\ \nonumber
& \Delta A\bar\Delta\delta(p_E)=\bar\Delta\delta(p_S) \\ \nonumber
& \bar\Delta\delta(p_S)=\bar\Delta\delta(p_W)\Delta B' \\
& \bar\Delta\delta(p_E)=\bar\Delta\delta(p_N)\Delta B \nonumber
\end{eqnarray}
If $w$ is the longest Weyl element, i.e. the anti-diagonal matrix with
only 1s on the anti-diagonal, then for a diagonal matrix $D$,
the matrix $wDw$ is again diagonal, with the sequence of the diagonal
elements reversed. The above equations then easily imply that
\[ \Delta A\Delta A'^{-1}=w \Delta B'\Delta B^{-1}w.\]
As these are all diagonal matrices this equation simply means the
sequence of equations
\[ \frac{a_j}{a'_j}=\frac{b'_{k+2-j}}{b_{k+2-j}},
\]
which is the periodicity property for the truncated $\Gamma$-system.

Finally, it easily follows from \eqref{xpera} that 
\[\frac{1-z_-}{1-z}=\frac{a'}a\quad\text{if } 1<j<k+1,
\]
 while according to \eqref{truncgam}, at the two extremes, we have
\[ \frac{1}{1-z}=\frac{a'}a\quad\text{if } j=1,\text{ and }
1-z_-=\frac{a'}a\quad\text{if } j=k+1.
\]

This immediately implies the periodicity for the $z$ variables lying
over the regular strip. As long as $r\geq k$, these ``regular''
variables clearly determine the remaining values, and thus the proof
is complete.

\vskip.4cm

\noindent Department of Geometry,
Institute of Mathematics, Budapest Institute of Technology, Budapest
H-1111, Hungary; {\tt szenes@math.bme.hu}

\end{document}